\documentclass{article}

\usepackage{
 amsmath,
 amsxtra,
 amsthm,
 amssymb,
 etex,
 mathrsfs,
 stmaryrd
 }
\usepackage[all]{xy}
\usepackage[T2A,T1]{fontenc}
\usepackage[OT2,T1]{fontenc}
\newcommand{\Zhe}{\mbox{\usefont{T2A}{\rmdefault}{m}{n}\CYRZH}}

\newtheorem{theorem}{Theorem}[section]
\newtheorem{lemma}[theorem]{Lemma}
\newtheorem{conjecture}[theorem]{Conjecture}
\newtheorem{proposition}[theorem]{Proposition}
\newtheorem{corollary}[theorem]{Corollary}
\newtheorem*{thm2}{Theorem}

\theoremstyle{definition}
\newtheorem{defn}[theorem]{Definition}
\newtheorem{remark}[theorem]{Remark}

\newcommand{\bd}{\begin{defn}}
\newcommand{\ed}{\end{defn}}
\newcommand{\bl}{\begin{lemma}}
\newcommand{\el}{\end{lemma}}
\newcommand{\bp}{\begin{proposition}}
\newcommand{\ep}{\end{proposition}}
\newcommand{\bt}{\begin{theorem}}
\newcommand{\et}{\end{theorem}}
\newcommand{\bc}{\begin{corollary}}
\newcommand{\ec}{\end{corollary}}
\newcommand{\br}{\begin{remark}}
\newcommand{\er}{\end{remark}}
\newcommand{\ba}{\begin{array}}
\newcommand{\ea}{\end{array}}
\newcommand{\bpf}{\begin{proof}}
\newcommand{\epf}{\end{proof}}

\newcommand{\Z}{\mathbb{Z}}
\newcommand{\Q}{\mathbb{Q}}
\newcommand{\Zp}{\mathbb{Z}_{p}}
\newcommand{\Qp}{\mathbb{Q}_{p}}
\newcommand{\Op}{\mathcal{O}}
\newcommand{\Ep}{E[p^{\infty}]}

\newcommand{\Ap}{A[p^{\infty}]}

\newcommand{\Ga}{\Gamma}

\newcommand{\La}{\Lambda}

\newcommand{\ord}{\mathrm{ord}}

\DeclareMathOperator{\Sel}{Sel} \DeclareMathOperator{\Gal}{Gal}
 \DeclareMathOperator{\rank}{rank}

\newcommand{\cyc}{\mathrm{cyc}}
\newcommand{\ac}{\mathrm{ac}}

\newcommand{\G}{\mathfrak{G}}

\newcommand{\mM}{\mathcal{M}}

\newcommand{\p}{\mathfrak{p}}
\newcommand{\q}{\mathfrak{q}}

\newcommand{\ot}{\otimes}
\newcommand{\ilim}{\displaystyle \mathop{\varinjlim}\limits}
\newcommand{\plim}{\displaystyle \mathop{\varprojlim}\limits}
\newcommand{\im}{\mathrm{im}\,}
\newcommand{\coker}{\mathrm{coker}\,}

\newcommand{\lra}{\longrightarrow}

\newcommand{\ps}[1]{\llbracket #1 \rrbracket}

  \DeclareFontFamily{U}{wncy}{}
  \DeclareFontShape{U}{wncy}{m}{n}{<->wncyr10}{}
  \DeclareSymbolFont{mcy}{U}{wncy}{m}{n}
  \DeclareMathSymbol{\sha}{\mathord}{mcy}{"58}

\numberwithin{equation}{section}

\begin{document}
\title{On fine Mordell-Weil groups over $\Zp$-extensions of an imaginary quadratic field}
 \author{Meng Fai Lim\footnote{School of Mathematics and Statistics, Central China Normal University, Wuhan, 430079, P.R.China. E-mail: \texttt{limmf@ccnu.edu.cn}} }
\date{}
\maketitle

\begin{abstract} \footnotesize
\noindent Let $E$ be an elliptic curve over $\Q$. Greenberg has posed a question whether the structure of the fine Selmer group over the cyclotomic $\Zp$-extension of $\Q$ can be described by cyclotomic polynomials in a certain precise manner. A recent work of Lei has made progress on this problem by proving that the fine Mordell-Weil group (in the sense of Wuthrich) does have this required property. The goal of this paper is to study analogous questions of Greenberg over various $\Zp$-extensions of an imaginary quadratic field $F$. In particular, when the elliptic curve has complex multiplication by the ring of integers of the imaginary quadratic field, we obtain analogous results of Lei over the cyclotomic $\Zp$-extension and anti-cyclotomic $\Zp$-extension of $F$. In the event that the elliptic curve has good ordinary reduction at the prime $p$, we further obtain a result over the $\Zp$-extension of $F$ unramified outside precisely one of the prime of $F$ above $p$. Finally, we study the situation of an elliptic curve over the anticyclotomic $\Zp$-extension under the generalized Heegner hypothesis. Along the way, we establish an analogous result for the BDP-Selmer group. This latter result is then applied to obtain a relation between the BDP $p$-adic $L$-function and the Mordell-Weil rank growth in the anticyclotomic $\Zp$-extension which may be of independent interest. 


\end{abstract}

\section{Introduction}

Throughout this article, $p$ will always denote a fixed prime $\geq 5$. The ($p$-primary)
fine Selmer group of an elliptic curve $E$ is a much studied object in Iwasawa theory
and frequently appears in the formulation (and proof) of the Iwasawa main conjecture (for instances, see \cite{CW, C83, K, Kob, PR00, Rub}). Despite this, it was only around 2005, Coates and Sujatha \cite{CS05}, and a little later in 2007, Wuthrich \cite{Wu} initiated a systematic study on the fine Selmer group of an elliptic curve $E$. Analogous to the usual Selmer group, the fine Selmer group $R(E/F)$ sits in the
middle of the following short exact sequence
\[0\lra \mM(E/F) \lra R(E/F) \lra \Zhe(E/F) \lra 0,\]
where $\mM(E/F)$ and $\Zhe(E/F)$ are the fine Mordell-Weil group and fine Tate-Shafarevich group respectively (defined in the sense of Wuthrich \cite{WuTS}) which can be thought as the ``fine'' counterpart of the usual Mordell-Weil group and Tate-Shafarevich group.

When the elliptic curve $E$ is defined over $\Q$, Greenberg asked whether the following formula
\[\mathrm{char}_{\Zp\ps{\Ga}} \Big(R(E/\Q^\cyc)^\vee\Big) = \prod_{n\geq 0} \Phi_n^{\max\{0, e_n-1\}}\]
is always valid (see \cite[Problem 0.7]{KurPo}). Here $e_n = \frac{\rank_\Z E(\Q_n)-\rank_\Z E(\Q_{n-1})}{\phi(p^n)}$ with $\phi(n)$ being the Euler totient function, $\Phi_n = \frac{(X+1)^{p^n}-1}{(X+1)^{p^{n-1}}-1}$ is the $p^n$th-cyclotomic polynomial, $\Q^\cyc$ is the cyclotomic $\Zp$-extension of $\Q$ and $\Q_n$ is the intermediate subextension of $\Q^\cyc/\Q$ with $|\Q_n:\Q| = p^n$. Pertaining to this question of Greenberg, Kurihara and Pollack (see loc.\ cit.) has performed several numerical computations verifying the truth of Greenberg's question. Furthermore, when the elliptic curve in question has good supersingular reduction at $p$, they have shown that the validity of Greenberg's question has implications towards understanding certain finer aspects of the structure of the signed Selmer groups of Kobayashi \cite{Kob} (also see \cite{LeiLim,LS} for some related discussion in this direction).

Recently, Lei \cite{LeiZ} has showed that
\[\mM(E/\Q^\cyc)^\vee \sim \prod_{n\geq 0} \Phi_n^{\max\{0, e_n-1\}}\]
under appropriate assumptions, where $\sim$ means pseudo-isomorphism of $\Zp\ps{\Gal(\Q^\cyc/\Q)}$-modules. Prior to this, Wuthrich has performed several numerical computations showing that $\Zhe(E/\Q^\cyc)$ is finite which led him to predict that this said group is finite most of the time (see \cite{WuTS}). Hence the collective works of Lei and Wuthrich provide strong evidence towards the validity of Greenberg's assertion. Inspired by these prior works, Lei went on to pose an analogous guess on the fine Selmer group of a $p$-supersingular CM elliptic curve over the anticyclotomic $\Zp$-extension in another paper of his \cite{LeiQ} but he did not provide an analogue of his above result in this context (we will see in a while a possible explanation for this). The goal of this paper is to continue this line of study. For the remainder of the introductory section, we describe our results in slightly more detail.

The theme of our work is whether one can formulate analogues of Greenberg's question for other different settings. The original work of Lei \cite{LeiZ} can be thought as providing evidences towards Greenberg's guess. Our approach towards this theme of study is to, a certain extent, the reversal of this development. Namely, we first seek to establish analogous results of Lei, and then propose forms of Greenberg's question suggested by the shape of our results.

A key principle towards our investigation is to understand the variation of the fine Mordell-Weil groups in the $\Zp$-extension of interest. In \cite{LeiZ}, Lei established a precise link between this said variation in the intermediate subextensions of the $\Zp$-extension $F_\infty/F$ and the $\Zp\ps{\Gal(F_\infty/F)}$-module structure of $\mM(E/F_\infty)$ under the assumption that the full Mordell-Weil group $E(F_\infty)$ over the $\Zp$-extension $F_\infty$ is finitely generated over $\Z$. However, this latter finite generation is not expected to always hold. For instance, if $E$ is an elliptic curve with complex multiplication given by the ring of integers of $F$ and has good supersingular reduction at $p$, then the Mordell-Weil group $E(F^\ac)$ is not finitely generated, where $F^\ac$ is the anticyclotomic $\Zp$-extension of an imaginary quadratic field $F$ (see \cite{AH2, CW, LLM}). In particular, the techniques developed in \cite{LeiZ} cannot apply in this context. In our paper, we circumvent this by instead establishing a relation between the structure of the fine Mordell-Weil groups in the intermediate subextensions of the $\Zp$-extension $F_\infty/F$ and the $\Zp\ps{\Gal(F_\infty/F)}$-module structure of an appropriate Iwasawa module $Y_f(E/F_\infty)$ (see Definition \ref{Yf def} and Proposition \ref{fineMWgeneral}) under the assumption that $R(E/F_\infty)^\vee$ is cotorsion over $\Zp\ps{\Gal(F_\infty/F)}$. This latter torsionness is conjectured to always hold (see \cite{LimFineDoc, LimFinePreprint, PR00}), and has been established in many cases of interest (see Remark \ref{torsion remark} and the references therein). We then show that our result recovers Lei's result (see Proposition \ref{fineMWgeneral2}) under his stronger assumption.

 Since the module $Y_f(E/F_\infty)$ is a quotient of $R(E/F_\infty)^\vee$ by definition, understanding this said module will provide information on the fine Selmer group. However, in view of the discussion in the preceding paragraph, the structure of $Y_f(E/F_\infty)$ will only become transparent after we are able to understand the explicit variation of the fine Mordell-Weil groups over the intermediate subextensions of the $\Zp$-extension in question.
The remainder of the paper is thus an effort towards realizing this latter aim. More precisely, we specialize to various $\Zp$-extensions of an imaginary quadratic field $F$, where we succeed in obtaining an explicit variation of the fine Mordell-Weil groups.
To illustrate, we present one of our results here. Suppose that $E$ is an elliptic curve with complex multiplication (which we sometimes abbreviate as ``CM elliptic curve'') given by the ring of integer of a imaginary quadratic field $F$ which has good reduction at $p$.
Suppose that $E$ is a CM elliptic curve defined over $\Q$ with given by $F$ and that the prime $p$ is inert in $F/\Q$. Write $\p$ for the unqiue prime of $F$ above $p$. Denote by $\Op$ (resp., $\Op_{\p})$ the ring of integers of $F$ (resp., the ring of integers of the local field $F_\p$). Let $F^\cyc$ be the cyclotomic $\Zp$-extension of $F$ and let $F_n$ denote the subextension of $F^\cyc/F$ such that $|F_n:F|=p^n$. Note that by our CM assumption, $E(F_n)$ comes equipped with a $\Op$-action, and so we can speak of $\rank_{\Op}E(F_n)$. Our first result is as follows. 

\begin{thm2}[Theorem \ref{mainCMcyc}]
  Retain settings as above. Write $\La_{\Op_\p}$ for $\Op_\p\ps{\Gal(F^\cyc/F)}$. Assume that $\Zhe(E/F_n)[p^\infty]$ is finite for every $n$. Then we have
\begin{equation}\label{mainCMcyceqn}
\mM(E/F^\cyc)^\vee\sim \bigoplus_{n\geq 0} (\La_{\Op_\p}/\Phi_n)^{\oplus \max\{0,f_n-1\}},
\end{equation}
where
\[ f_n = \left\{
           \begin{array}{ll}
            \displaystyle\frac{\rank_{\Op}E(F_n) - \rank_{\Op}E(F_{n-1}) }{\phi(p^n)} , & \mbox{if $n\geq 1$;} \\
            \\
            \rank_{\Op}E(F) , & \mbox{if $n=0$.}
           \end{array}
         \right.
\]
\end{thm2}
In this supersingular reduction situation, we also have a result for the anticyclotomic $\Zp$-extension (see Theorem \ref{mainCManticyc}). In the event that the CM elliptic curve has good ordinary reduction at $p$, besides the cyclotomic $\Zp$-extension (see Theorem \ref{mainCMpsplitcyc}) and anticyclotomic $\Zp$-extension (see Theorem \ref{mainCMpsplitac}), we also obtain a result over the $\Zp$-extension unramified outside a fixed prime of $F$ above $p$ (cf.\ Theorem \ref{mainCMpsplit1}). Motivated by these results, we posed analogous questions of Greenberg in each of these contexts (see Questions 1, 2, 3, 4 and 5).

We finally consider the situation of an elliptic curve with good reduction at $p$ over the anticyclotomic $\Zp$-extension of an imaginary quadratic field $F$ at which $p$ splits completely and a certain generalized Heegner condition holds (see first paragraph of Section \ref{BDP section}). In this context, our result on the fine Mordell-Weil groups (see Theorem \ref{BDPfine}) is slightly less precise but it's enough to inspire us to pose an analogous question of Greenberg in this context (see Question 6). Along the way, we will prove an analogous result for the BDP-Selmer group (see Theorem \ref{BDPac}). This latter theorem, although not directly related to fine Selmer groups, seems to be interesting in its own right and is established using the methods of this paper. This theorem also has an application towards yielding a relation between the BDP $p$-adic $L$-function and the Mordell-Weil rank growth in the anticyclotomic $\Zp$-extension (see Corollary \ref{BDPacCor}).

We now give an outline of the paper. In Section \ref{Preliminaries}, we first collect certain facts from \cite{Lee} which we then apply to obtain a general algebraic result in the form of Corollary \ref{tech module}. This said result will provide the necessary foundation for us to obtain our Proposition \ref{fineMWgeneral} which links the structure of the fine Mordell-Weil groups at intermediate subextensions of the $\Zp$-extension $F_\infty/F$ to the $\Zp\ps{\Gal(F_\infty/F)}$-module structure of the Iwasawa module $Y_f(E/F_\infty)$. We also show how a result of Lei can be recovered from ours. In Section \ref{CMinert section}, we specialize to the context of a CM elliptic curve at which the prime $p$ is a good supersingular one. In particular, we are able to provide a theoretical evidence towards Lei's guess in \cite{LeiQ}. We then turn to the context of a CM elliptic curve at which the prime $p$ is a good ordinary prime in Section \ref{CMsplit section}. Finally, in Section \ref{BDP section}, we consider the context of an elliptic curve with good reduction at $p$ under the generalized Heegner hypothesis.

\subsection*{Acknowledgement}
The author would like to thank Antonio Lei for many insightful
discussions and comments on the paper. He would also like to thank the anonymous referee for the many suggestions and comments on the paper. 
  The research of the author is supported by the
National Natural Science Foundation of China under Grant No. 11771164.

\section{Algebraic preliminaries} \label{Preliminaries}

In this section, we recall certain algebraic preliminaries. We begin by introducing certain basic conventions. For a $\Zp$-module $N$, we denote by $N[p^n]$ the submodule of $N$ consisting of elements
of $N$ which are annihilated by $p^n$. We then write $N[p^\infty] = \cup_{n\geq 1}N[p^n]$ and $T_pN = \plim_n N[p^n]$. In the event that $N$ is a finite abelian $p$-group, we shall write $\ord_p(N)$ for the power of $p$ in $|N|$, i.e., $|N| = p^{\ord_p(N)}$. Finally, for an abelian group or module $M$, we often write $M^{\oplus n}$ for the $n$-copies of $M$.

\subsection{The functor $\mathfrak{G}$}

Throughout the paper, we shall write $\La$ for the classical Iwasawa algebra $\Zp\ps{\Ga}$, where $\Ga=\Zp$. Recall that a homomorphism $M\lra N$ of $\La$-modules is said to be a 
pseudo-isomorphism if its kernel and cokernel are finite. We shall write $M\stackrel{\sim}{\lra} N$. Note that the notion of being pseudo-isomorphic is not a reflexive relation (see \cite[P 300, Exercise 1]{NSW}). However, in the event that $M$ and $N$ are torsion $\La$-modules, then one has $M\stackrel{\sim}{\lra} N$ if and only if $N\stackrel{\sim}{\lra} M$ (see \cite[Remarks after Proposition 5.1.7]{NSW}). In this case, we shall sometimes write $M\sim N$.

The unique subgroup of $\Ga$ of index $p^n$ is in turn denoted by $\Ga_n$. For a given $\La$-module $M$, Lee defined a corresponding $\La$-module
\[\G(M): = \plim_n \Big(M_{\Ga_n}[p^\infty]\Big). \]
For more details on this module, the interested readers are referred to the paper of Lee \cite{Lee}.
For our purposes, we require the following.

\bl \label{G functor}
Let $M$ be a finitely generated $\La$-module. Then $\G(M)=0$ if and only if there is a pseudo-isomorphism
  \[M\stackrel{\sim}{\lra}\La^r \oplus \Big(\bigoplus_{j=1}^s \La/\Phi_{m_j}\Big),\]
  where $\Phi_{m_j}$ is a certain $p^{m_j}$-th cyclotomic polynomial.
\el

\bpf
See \cite[Lemma 4.1.3]{Lee}.
\epf

\subsection{Some technical results}

We now develop an axiomatic setting which will be useful for our eventual arithmetic application. As before, we let $\La$ denote the classical Iwasawa algebra $\Zp\ps{\Ga}$, where $\Ga\cong \Zp$. Let $\{M_n\}$ be a sequence of $\La$-modules with transition maps $M_{n+1}\lra M_n$ such that the action of $\La$ on $M_n$ factors through $\Zp[\Ga/\Ga_n]$. We shall always assume that each $M_n$ is finitely generated over $\Zp$.

Writing $M_{n,f} = M_n/M_n[p^\infty]$, one can check easily that every $M_n[p^\infty]$ (resp., $M_{n,f}$) is a $\La$-submodule (resp., $\La$-quotient) of $M_n$ and the transition map $M_{n+1}\lra M_n$ induces a transition map $M_{n+1}[p^\infty]\lra M_n[p^\infty]$ (resp., $M_{n+1,f}\lra M_{n,f}$).
Set $N= \plim_n (M_n[p^\infty])$, $M = \plim_n M_n$ and $M_f = \plim_n M_{n,f}$.
For each $n$, the natural map $M\lra M_n$ factors through $M_{\Ga_n}$ to induce a map $M_{\Ga_n}\lra M_n$. Similarly, we have maps $N_{\Ga_n}\lra M[p^\infty]$ and $(M_f)_{\Ga_n}\lra M_{n,f}$.

\br \label{warning1}
We warn the readers that  $\plim_n (M_n[p^\infty])$ and  $(\plim_n M_n)[p^\infty]$ need not be equal and so $M_f$ needs not be the same as $M/M[p^\infty]$. We shall illustrate this with a class of examples. Set $M=\La/g(T)$, where $g(T)$ is any irreducible polynomial coprime to every $p$-power cyclotomic polynomials. Then set $M_n = M_{\Ga_n}$. Plainly, one has $M=\plim_n M_n$, $M_n = M_n[p^\infty]$ and $M_{n,f} =0$. Thus, by definition, we have $M_f =0$. On other hand, $M=\La/g(T)$ has no $p$-torsion, and so one has $M[p^\infty] =0$ and $M/M[p^\infty] = M \neq 0$.
\er

\bp \label{G=0}
Retain settings as above. Suppose further that the map $M_{\Ga_n}\lra M_n$ has finite kernel and cokernel for every $n$. Then one has $\G(M_f)=0$. In other words, one has
\[ M_f \stackrel{\sim}{\lra} \La^{\oplus r} \oplus \Big(\bigoplus_{n\geq 0} (\La/\Phi_n)^{\oplus s_n}\Big),\]
where $r$ is a nonnegative integer and $\{s_n\}$ is a sequence of nonnegative integers with $s_n=0$ for $n\gg 0$.
\ep

\bpf
Consider the following commutative diagram
\[   \entrymodifiers={!! <0pt, .8ex>+} \SelectTips{eu}{}\xymatrix{
     & N_{\Ga_n} \ar[d] \ar[r] &  M_{\Ga_n}
    \ar[d] \ar[r] & (M_f)_{\Ga_n} \ar[d]^{g_n} \ar[r] & 0\\
    0 \ar[r]^{} & M_n[p^\infty] \ar[r]^{} & M_n \ar[r] & M_{n,f} \ar[r] & 0
     } \]
with exact rows. In view of our standing assumption of $M_n$ being finitely generated over $\Zp$, the module $M_n[p^\infty]$ is therefore finite. It then follows from this and the hypothesis of the proposition that $\ker g_n$ is finite for every $n$. Since $M_{n,f}$ is $\Zp$-torsionfree by definition, this in turn implies that
\[ (M_f)_{\Ga_n}[p^\infty] = \ker g_n. \]
Taking inverse limit, we obtain
\[ \G(M_f) \cong \plim_n \ker g_n\]
but the latter is trivial by virtue of $\plim_n (M_f)_{\Ga_n} \cong M_f = \plim_n M_{n,f}$. This establishes the first assertion of the proposition. The second is now immediate from this and Lemma \ref{G functor}.
\epf

The following corollary will be a crucial component for our subsequent discussion in the paper.

\bc \label{tech module}
Retain the settings and assumptions in Proposition \ref{G=0}. Suppose further that the following conditions are valid.
\begin{enumerate}
  \item[$(a)$] The module $M$ is torsion over $\La$. In particular, one has a pseudo-isomorphism
  \[ M_f\sim \bigoplus_{n\geq 0} (\La/\Phi_n)^{\oplus s_n}\]
of $\La$-modules, where $\{s_n\}$ is a sequence of nonnegative integers with $s_n=0$ for $n\gg 0$.
  \item[$(b)$] $\{M'_n\}$ is a sequence of $\La$-modules with transition maps $M'_{n+1}\lra M'_n$ such that the action of $\La$ on $M'_n$ factors through $\Zp[\Ga/\Ga_n]$ and that each $M'_n$ is a $\La$-submodule of $M_n$ with finite cardinality.
\end{enumerate}
 Then for every $n$, one has
 \[ T_p\big((M_n/M'_n)^\vee\big) \cong \bigoplus_{j= 0}^n (\La/\Phi_j)^{\oplus s_j},\]
 where the integer $n$ appearing on both sides of the isomorphism are the same. \ec

\bpf
Since $M'_n\subseteq M_n[p^\infty]$, we have the following commutative diagram
\[   \entrymodifiers={!! <0pt, .8ex>+} \SelectTips{eu}{}\xymatrix{
    0 \ar[r] & M'_n \ar[d] \ar[r] &  M_n
    \ar@{=}[d] \ar[r] & M_n/M'_n \ar[d]\ar[r] & 0\\
    0 \ar[r]^{} & M_n[p^\infty] \ar[r]^{} & M_n \ar[r] & M_{n,f} \ar[r] & 0
     } \]
with exact rows. By the snake lemma, one has an exact sequence
\[ 0\lra M'_n\lra M_n[p^\infty] \lra M_n/M'_n \lra M_{n,f}\lra 0.\]
Since $M_n[p^\infty]$ is finite and the Tate module of a finite group is trivial, we therefore have $T_p\big((M_n/M'_n)^\vee\big) \cong  T_p\big(M_{n,f}^\vee\big)$. On the other hand, it follows from the proof of the preceding proposition that the map
$(M_f)_{\Ga_n} \lra M_{n,f}$
has finite kernel and cokernel. This in turn yields $T_p\big(M_{n,f}^\vee\big) \cong T_p\big(((M_f)^\vee)_{\Ga_n}\big)$. By the structural result on $M_f$ in the preceding proposition, one then has $T_p\big(((M_f)^\vee)_{\Ga_n}\big) \cong \bigoplus_{j= 0}^n (\La/\Phi_j)^{\oplus s_j}$. Thus, we have the conclusion of the corollary.
\epf

\section{Fine Selmer groups over $\Zp$-extension}

We begin with a remark. Although the paper is concerned with elliptic curves, the discussion in this section also holds for abelian varieties, and whence our presentation in this section will take this general approach. Now, let $A$ be an abelian variety defined over a number field $F$. Let $S$ be a finite set of primes of $F$ containing the primes above $p$, the bad reduction primes of $A$ and the infinite primes. Denote by $F_S$ the maximal algebraic extension of $F$ which is unramified outside $S$. For every extension $\mathcal{L}$ of $F$ contained in $F_S$, write $G_S(\mathcal{L})=\Gal(F_S/\mathcal{L})$, and denote by $S(\mathcal{L})$ the set of primes of $\mathcal{L}$ above $S$.
We shall also write $S_p(\mathcal{L})$ for the set of primes of $\mathcal{L}$ above $p$.

Let $L$ be a finite extension of $F$ contained in $F_S$. Following Coates, Sujatha and Wuthrich \cite{CS05, Wu}, the fine Selmer group $R(A/L)$ of $A$ over $L$ is defined by the exact sequence
\[0\lra R(A/L)\lra H^1(G_S(L),\Ap)\lra \bigoplus_{v\in S(L)} H^1(L_v, \Ap).\]

At first viewing, it would seem that the fine Selmer group depends on the set $S$. But we shall see that this is not so. In fact, recall that the classical ($p$-primary) Selmer group $\Sel(A/L)$ is defined by the exact sequence
\[0\lra \Sel(A/L)\lra H^1(G_S(L),\Ap)\lra \bigoplus_{v\in S(L)} H^1(L_v, A)[p^\infty]\]
and it is well-known that this definition is independent of the set $S$ as long as the set $S$ contains all the primes above $p$ and the bad reduction primes of $A$ (see \cite[Chap. I, Corollary 6.6]{Mi}). Furthermore, we have a short exact sequence
\begin{equation} \label{eqn short exact} 0 \lra A(L)\ot_{\Zp}\Qp/\Zp \lra \Sel(A/L)\lra \sha(A/L)[p^\infty]\lra 0, \end{equation}
where $\sha(A/L)$ is the Tate-Shafarevich group.

 The fine Selmer group and the classical Selmer group are related by the following exact sequence.

\bl
We have an exact sequence
\[ 0\lra R(A/L) \lra \Sel(A/L) \lra \bigoplus_{v\in S_p(L)}A(L_v)\ot_{\Zp}\Qp/\Zp.\]
In particular, the definition of the fine Selmer group does not depend on the
choice of the set $S$. \el

\bpf
 See \cite[Lemma 4.1]{LMu}.
\epf

We now define the fine Mordell-Weil group and fine Tate-Shafarevich group following Wuthrich \cite{WuTS}. As a start, the fine Mordell-Weil group $\mathcal{M}(A/L)$ is defined by
\[ \mathcal{M}(A/L) = \ker\Big(A(L)\ot_{\Zp}\Qp/\Zp \lra \bigoplus_{v\in S_p(L)} A(L_v)\ot_{\Zp}\Qp/\Zp \Big),\]
which fits into the following commutative diagram
\[   \entrymodifiers={!! <0pt, .8ex>+} \SelectTips{eu}{}\xymatrix{
    0 \ar[r]^{} & \mM(A/L) \ar[d] \ar[r] &  A(L)\ot_{\Zp}\Qp/\Zp
    \ar[d] \ar[r] & \displaystyle\bigoplus_{v\in S_p(L)} A(L_v)\ot_{\Zp}\Qp/\Zp  \ar@{=}[d]\\
    0 \ar[r]^{} & R(A/L) \ar[r]^{} & \Sel(A/L)\ar[r] & \displaystyle\bigoplus_{v\in S_p(L)} A(L_v)\ot_{\Zp}\Qp/\Zp
     } \]
with exact rows, where the leftmost vertical map is induced by the middle vertical map. The fine Tate-Shafarevich group $\Zhe(A/L)$ is then defined by
\[ \Zhe(A/L) = \coker\Big( \mM(A/L)\lra R(A/L)\Big).\]
Since the middle vertical map in the above diagram is injective, so is the leftmost vertical map. Consequently, via a snake lemma argument, we have a short exact sequence
\begin{equation} \label{eqn fine short exact} 0 \lra \mM(A/L) \lra R(A/L) \lra \Zhe(A/L)\lra 0
\end{equation}
with $\Zhe(A/L)$ injecting into $\sha(A/L)[p^\infty]$.

For a given $\Zp$-extension $F_\infty$ of $F$, we denote by $\Ga$ the Galois group of the extension $F_\infty/F$ and let $\La$ denote the resulting Iwasawa algebra over $\Zp$. We in turn denote by $F_n$ the fixed field of $\Ga_n$, where $\Ga_n$ is the unique subgroup of $\Ga$ of index $p^n$. The fine Selmer group of $A$ over $F_\infty$ is then defined to be $R(A/F_\infty) = \ilim_n R(A/F_n)$ which comes naturally equipped with a $\La$-module. The $\La$-modules $\mM(A/F_\infty)$ and $\Zhe(A/F_\infty)$ are defined via similar limiting processes. We then write $Y(A/F_n)$ and $Y(A/F_\infty)$ for the Pontryagin dual of $R(A/F_n)$ and $R(A/F_\infty)$ respectively.

To facilitate subsequent discussion, we recall the control theorem for fine Selmer group (see \cite[Theorem 3.3]{LimFineDoc} for the proof).

\bt \label{control theorem}
Let $A$ be an abelian variety defined over $F$ and let $F_{\infty}$ be a $\Zp$-extension of $F$.
Then the restriction map
\[r_n: R(A/F_n) \lra R(A/F_{\infty})^{\Ga_n}\]
has finite kernel and cokernel which is bounded independent of $n$.
\et

We also recall the following conjecture (see \cite{LimFineDoc, LimFinePreprint, PR00}).

 \medskip
 \noindent
\textbf{Conjecture Y.} Let $A$ be an abelian variety defined over a number field $F$ and $F_\infty$ a $\Zp$-extension of $F$. Then $Y(A/F_\infty)$ is torsion over $\Zp\ps{\Gal(F_\infty/F)}$.

\br \label{torsion remark}
In the event that $E$ is an elliptic curve over $\Q$ with good reduction at $p$ and $F$ is an abelian extension of $\Q$, then Conjecture Y is valid for $Y(E/F^\cyc)$ by a theorem of Kato \cite{K}. There are also some previously known cases for non-cyclotomic $\Zp$-extensions (see \cite{AH2, AH, C83, PR00}). We also refer readers to \cite[Corollary 3.5 and Remark 5.2]{LimFineDoc} and \cite[Theorem 3.9 and Section 6]{LimFinePreprint} for more theoretical and numerical evidences on this conjecture.
\er

\bd \label{Yf def} Set $Y_f(A/F_\infty) = \plim_n Y_f(A/F_n)$, where $Y_f(A/F_n) = Y(A/F_n)/ Y(A/F_n)[p^\infty]$. \ed

\br
We warn the readers that $Y_f(A/F_\infty)$ may not coincide with $Y(A/F_\infty)/Y(A/F_\infty)[p^\infty]$ in general (see Remark \ref{warning1}).
\er

Now, supposing that $Y(A/F_\infty)$ is torsion over $\Zp\ps{\Gal(F_\infty/F)}$, it then follows from a combination of Theorems \ref{G=0} and \ref{control theorem} that
\begin{equation}\label{strucDiv}
  Y_f(A/F_\infty)\sim \bigoplus_{n\geq 0} (\La/\Phi_n)^{\oplus s_n},
\end{equation}
where $\{s_n\}$ is a sequence of nonnegative integers with $s_n = 0$ for $n\gg 0$.

\bp \label{fineMWgeneral}
Let $A$ be an abelian variety defined over $F$ and let $F_{\infty}$ be a $\Zp$-extension of $F$.
Suppose that the following statements are valid.
\begin{enumerate}
  \item[$(a)$] $\Zhe(A/F_n)$ is finite for every $n$.
  \item[$(b)$] $Y(A/F_\infty)$ is torsion over $\Zp\ps{\Gal(F_\infty/F)}$.
\end{enumerate}
Then we have
\[T_p(\mM(A/F_n)) \cong \bigoplus_{j= 0}^n (\La/\Phi_j)^{\oplus s_j}\]
for every $n$, where the $s_j$'s are given as in $(\ref{strucDiv})$.
\ep

\bpf
Set $M_n = Y(A/F_n)$ and $M'_n = \Zhe(A/F_n)^\vee$. Therefore, $M_n/M'_n = \mM(A/F_n)^\vee$. By definition, $M = \plim_n M_n = Y(A/F_\infty)$. Now note that the map $M_{\Ga_n}\lra M_n$ has finite kernel and cokernel by Theorem \ref{control theorem}. Therefore, we may apply Corollary \ref{tech module} to obtain the conclusion of the proposition.
\epf

A variant of Proposition \ref{fineMWgeneral} was obtained by Lei for the fine Mordell-Weil group in \cite[Corollary 3.8]{LeiZ} under the assumption that $A(F_\infty)$ is finitely generated over $\Z$. We now show how our result may recover Lei's result under his assumption.

\bp \label{fineMWgeneral2}
Let $A$ be an abelian variety defined over $F$ and let $F_{\infty}$ be a $\Zp$-extension of $F$.
Suppose that the following statements are valid.
\begin{enumerate}
  \item[$(a)$] $\Zhe(A/F_n)$ is finite for every $n$.
  \item[$(b)$] $A(F_\infty)$ is finitely generated over $\Z$.
\end{enumerate}
Then we have
\[ \mM(A/F_\infty)^\vee\sim \bigoplus_{n\geq 0} (\La/\Phi_n)^{\oplus s_n}\]
with the property that
\[T_p(\mM(A/F_n)) \cong \bigoplus_{j\geq 0}^n (\La/\Phi_j)^{\oplus s_j}\]
for every $n$.
\ep

\bpf
Under the assumption that $A(F_\infty)$ is finitely generated over $\Z$, it follows from a result of Lee \cite[Theorem 2.1.2]{Lee} that
\begin{equation}\label{mM}
  \mM(A/F_\infty)^\vee\sim \bigoplus_{n\geq 0} (\La/\Phi_n)^{\oplus r_n},
\end{equation}
where $\{r_n\}$ is a sequence of nonnegative integers with $r_n = 0$ for $n\gg 0$. On the other hand, it follows from \cite[Corollary 3.5]{LimFinePreprint} that the hypotheses of the proposition imply that $Y(A/F_\infty)$ is torsion over $\La$, and so Proposition \ref{fineMWgeneral} applies. It therefore remains to show that $\mM(A/F_\infty)^\vee \sim Y_f(A/F_\infty)$. Now for each $n$, it follows from the argument in the proof of Corollary \ref{tech module} that there is an exact sequence
\[ 0\lra \Zhe(A/F_n)^\vee \lra Y(A/F_n)[p^\infty] \lra \mM(A/F_n)^\vee \stackrel{h_n}{\lra} Y_f(A/F_n)\lra 0 \]
for every $n$. In particular, it follows that $\ker h_n$ is finite and so is contained in $\mM(A/F_n)^\vee[p^\infty]$. As a consequence, we have
\[ \ord_p\big(\ker h_n \big) \leq \ord_p\Big(\mM(A/F_n)^\vee[p^\infty]\Big) \leq \ord_p\Big(\big(\mM(A/F_\infty)^\vee \big)_{\Ga_n}[p^\infty]\Big) +O(1),\]
where the second inequality follows from the observation that the kernel of the map
\[ \mM(A/F_n) \lra  \mM(A/F_\infty)^{\Ga_n} \]
is contained in the kernel of the map
\[ R(A/F_n) \lra  R(A/F_\infty)^{\Ga_n} \]
which in turn is finite and bounded independently of $n$ by Theorem \ref{control theorem}. Now, note that by Lemma \ref{G functor} and (\ref{mM}), one has $\G\big(\mM(A/F_\infty)^\vee\big)=0$. It then follows from this and an application of \cite[Lemma 4.1.3]{Lee} that
\[ \ord_p\Big(\big(\mM(A/F_\infty)^\vee \big)_{\Ga_n}[p^\infty]\Big) = O(1).\]
Hence we may conclude that $\ker h_n$ is finite and bounded independently of $n$. Consequently, upon taking inverse limit, we see that the map
\[ \mM(A/F_\infty)^\vee \lra Y_f(A/F_\infty)\]
is surjective with a finite kernel which establishes the desired pseudo-isomorphism.
\epf

\subsection{A remark on Mordell-Weil growth in $\Zp$-extension}

In this subsection, we shall study the variation of Mordell-Weil groups in a $\Zp$-extension. Our discussion here will rely heavily on a result of Lee \cite{Lee}. The result we seek to prove is as follow.

\bp \label{MWgrowth}
Let $A$ be an abelian variety defined over a number field $F$ and $F_\infty/F$ a $\Zp$-extension. Define
\[ e_n = \left\{
           \begin{array}{ll}
            \displaystyle\frac{\rank_{\Z}A(F_n) - \rank_{\Z}A(F_{n-1}) }{\phi(p^n)} , & \mbox{if $n\geq 1;$} \\
            \\
            \rank_{\Z}A(F) , & \mbox{if $n=0$.}
           \end{array}
         \right.
\]
Then we have
\[ T_p\big(A(F_n)\ot \Qp/\Zp \big) \cong \bigoplus_{j=0}^n(\La/\Phi_j)^{\oplus e_j}.\]
\ep

This proposition might perhaps be well-known among the experts but due to a lack of reference, we shall supply a proof here. Note that the proposition \textit{does not} have any assumption on the reduction type of the abelian variety $A$ or the $\Zp$-extension in question. We also emphasize that the proposition \textit{does not} assume the finiteness of $\sha(A/F_n)[p^\infty]$. 

\bpf[Proof of Proposition \ref{MWgrowth}]
As a start, we note that all intermediate objects $A(F_n)\ot\Qp/\Zp$, $(A(F_\infty)\ot\Qp/\Zp)^{\Ga_n}$ and their Pontryagin dual will be viewed as $\La$-modules without further mention.
By a theorem of Lee \cite[Theorem 2.1.2]{Lee}, there is an injective $\La$-homomorphism
\[ \big(A(F_\infty)\ot\Qp/\Zp\big)^\vee \lra \La^{\oplus r}\oplus \Big(\bigoplus_{j\geq 0}(\La/\Phi_j)^{\oplus t_j}\Big)\]
with a finite cokernel,
where $\{t_j\}$ is a sequence of nonnegative integers with $t_j = 0$ for $j\gg0$.
An application of the $-_{\Ga_n}$-functor yields a $\La$-map
\begin{equation}\label{A MW1}
\Big((A(F_\infty)\ot\Qp/\Zp)^{\Ga_n}\Big)^\vee \lra \bigoplus_{j\geq 0}(\La/\Phi_j)^{\oplus (r+ t_j)}
\end{equation}
with finite kernel and cokernel. On the other hand, the restriction map
\[ A(F_n)\ot\Qp/\Zp \lra (A(F_\infty)\ot\Qp/\Zp)^{\Ga_n}\]
has a kernel contained in $H^1(\Ga_n, A(F_\infty)[p^\infty])$, where the latter is finite by \cite[Lemma 3.4]{LimFineDoc}. Consequently, we have a $\La$-map
\begin{equation}\label{A MW2} \Big((A(F_\infty)\ot\Qp/\Zp)^{\Ga_n}\Big)^\vee \lra \big(A(F_n)\ot\Qp/\Zp\big)^\vee.
\end{equation}
with finite cokernel.  Now, a combination of (\ref{A MW1}) and (\ref{A MW2}) will yield
\[ \big(A(F_n)\ot\Qp/\Zp\big)^\vee\sim \bigoplus_{j=0}^n(\La/\Phi_j)^{\oplus c^{(n)}_j} \]
for some nonnegative integers $c^{(n)}_j$, where the integers may possibly depend on $n$. (The point is that (\ref{A MW1}) and (\ref{A MW2}) do not a priori tell us whether the number of occurrences of a particular $\Phi_j$ in $\big(A(F_m)\ot\Qp/\Zp\big)^\vee$ is the same as that in $\big(A(F_n)\ot\Qp/\Zp\big)^\vee$ for sufficiently large distinct $m$ and $n$.)

We now set to show that the $c_j^{(n)}$'s do not depend on $n$. Suppose that $m>n$. Consider the restriction map
\[ A(F_n)\ot\Qp/\Zp \lra \big(A(F_m)\ot\Qp/\Zp\big)^{\Ga_n} = \big(A(F_m)\ot\Qp/\Zp\big)^{\Ga_n/\Ga_m}\]
which fits into the following commutative diagram
\[   \entrymodifiers={!! <0pt, .8ex>+} \SelectTips{eu}{}\xymatrix{
    0 \ar[r]^{} & A(F_n)\ot\Qp/\Zp \ar[d] \ar[r] &  H^1(F_n, \Ap)
    \ar[d]^{h_{m,n}} \ar[r] & H^1(F_n, A)[p^\infty] \ar[d]^{g_{m,n}} \ar[r] & 0\\
    0 \ar[r]^{} & \big(A(F_m)\ot\Qp/\Zp\big)^{\Ga_n/\Ga_m} \ar[r]^{} & H^1(F_m, \Ap)^{\Ga_n/\Ga_m}\ar[r] & H^1(F_m, A)[p^\infty]^{\Ga_n/\Ga_m} &
     } \]
      Via a Hochschild-Serre spectral sequence argument, we see that $\ker h_{m,n} = H^1(\Ga_n/\Ga_m, A(F_m)[p^\infty])$, $\ker g_{m,n} = H^1(\Ga_n/\Ga_m, A(F_m))$   and $\mathrm{coker}~h_{m,n} \subseteq H^2(\Ga_n/\Ga_m, A(F_m)[p^\infty])$. Since $\Ga_n/\Ga_m$ and $A(F_m)[p^\infty]$ are finite, so are $\ker h_{m,n}$ and $\mathrm{coker} ~h_{m,n}$. On the other hand, as $A(F_m)$ is a finitely generated abelian group by Mordell-Weil Theorem, we see that $\ker g_{m,n} = H^1(\Ga_n/\Ga_m, A(F_m))$ is also finite.
Hence it follows that the kernel and cokernel of the leftmost vertical map are finite. Consequently, we have
\[\Big(A(F_n)\ot\Qp/\Zp\Big)^\vee \sim \Big(\big((A(F_m)\ot\Qp/\Zp)\big)^\vee\Big)_{\Ga_n},\]
or equivalently
\[ \bigoplus_{j=0}^n(\La/\Phi_j)^{\oplus c^{(n)}_j} \sim \bigoplus_{j=0}^n(\La/\Phi_j)^{\oplus c^{(m)}_j} .\]
This in turn shows that $c^{(n)}_j= c^{(m)}_j$ for $0\leq j\leq n$. Hence, in conclusion, we have shown that there exists a sequence $\{c_j\}$ of nonnegative integers such that
\[ \big(A(F_n)\ot\Qp/\Zp\big)^\vee\sim \bigoplus_{j=0}^n(\La/\Phi_j)^{\oplus c_j}. \]
for every $n$. Now, by an iterative rank comparison, it follows that $c_n =e_n$ for every $n$. On the other hand, observe that the Tate module of a finite group is trivial and the Tate module of a finitely generated $\Zp$-free module is itself. Hence, by combining this latter observation with the above pseudo-isomorphism, we obtain the isomorphism as asserted in Proposition \ref{MWgrowth}.
\epf

We also state the following local analog, where a proof is again supplied for the convenience of the readers.

\bp \label{MWgrowthlocal}
Let $K$ be a finite extension of $\Qp$ and $K_\infty$ a $\Zp$-extension of $K$.
Let $A$ be a $g$-dimensional abelian variety defined over $K$.
 Then we have
\[ T_p\big(A(K_n)\ot \Qp/\Zp \big) \cong \bigoplus_{j=0}^n(\La/\Phi_j)^{\oplus g|K:\Q_p|}.\]
\ep

\bpf
By Lee's theorem \cite[Theorem 2.1.2]{Lee}, one has
\[ \big(A(K_\infty)\ot\Qp/\Zp\big)^\vee \stackrel{\sim}{\lra} \La^{\oplus r}\oplus \Big(\bigoplus_{j\geq 0}(\La/\Phi_j)^{\oplus t_j}\Big).\]
Therefore, a similar argument to that in Proposition \ref{MWgrowth} (where we replace the application of Mordell-Weil Theorem with Mattuck's theorem \cite{Mat}) shows that there is a sequence $\{d_j\}$ of nonnegative integers such that
\[ \big(A(K_n)\ot\Qp/\Zp\big)^\vee\sim \bigoplus_{j=0}^n(\La/\Phi_j)^{\oplus d_j}. \]
But Mattuck's theorem also tells us that $\mathrm{rank}_{\Zp}A(K_n) = g|K_n:\Qp| = gp^n|K:\Qp|$. Taking this into account, it then follows by an induction argument that $d_n = g|K:\Qp|$ for every $n$. The conclusion of the proposition is now immediate from this. \epf

\section{CM elliptic curve: good supersingular reduction case} \label{CMinert section}

From now on, $F$ will always denote an imaginary quadratic field with ring of integers $\Op:=\Op_F$. Let $E$ be an elliptic curve defined over $\Q$ with complex multiplication given by $\Op$.  Note that $F$ necessarily has class number $1$. As before, we let $p$ be a prime $\geq 5$. \textbf{In this section, the elliptic curve $E$ is always assumed to have good supersingular reduction at $p$, and so the prime $p$ is inert in $F/\Q$.}  We then denote by $\p$ the unique prime of $F$ above $p$ and write $\Op_\p:=\Op_{F_\p}$ for the ring of integers of the local field $F_\p$. Note that in this context, $\Ep$ comes naturally endowed with a $\Op_\p[\Gal(\overline{\Q}/F)]$-module structure. Moreover, as an $\Op_\p$-module, it is cofree of corank $1$. Consequently, for every finite extension $L$ of $F$, the fine Selmer group $R(E/L)$ has a natural $\Op_\p$-module structure. In particular, if $F_\infty$ is a $\Zp$-extension of $F$ with Galois group $\Ga$, $R(E/F_\infty)$ has a $\La_{\Op_\p}$-module structure, where $\La_{\Op_\p}:=\Op_\p\ps{\Ga}$. We have similar assertions for the fine Mordell-Weil group and fine Tate-Shafarevich group.

On the other hand, if $L$ is a finite Galois extension of our imaginary quadratic field $F$, the Mordell-Weil group $E(L)$ has a natural $\Op$-module structure. Therefore, it makes sense to speak of $\rank_{\Op}E(L)$ and one clearly has $\rank_{\Op}E(L) = 2\rank_{\Z}E(L)$. Furthermore, we can speak of $E(L)\ot_{\Op} F_\p/\Op_\p$. Now, since $p$ is inert in $F/\Q$, one has $\p = p\Op$ and so the Tate module of $E(L)\ot_{\Op} F_\p/\Op_\p$ with respect to $\p$ coincides with the usual Tate module of $E(L)\ot_{\Op} F_\p/\Op_\p$ with respect to $p$. We shall therefore write $T_p\big(E(L)\ot_{\Op} F_\p/\Op_\p\big)$ for this Tate module which comes equipped with an $\Op_\p$-action. As a consequence, we can speak of $T_p\big(E(L)\ot_{\Op} F_\p/\Op_\p\big) \ot_{\Op_\p}F_\p$.

Now, assuming further that $\p$ is totally ramified in $L/F$, we shall, by abuse of notation, write $\p$ for the prime of $L$ above $\p$ (and hence $p$). Then there is a natural inclusion
$E(L) \hookrightarrow E(L_{\p})$ 
 which in turn induces a $F_\p[\Gal(L/F)]$-map
\begin{equation}\label{f_L}
 f_L : T_p\big(E(L)\ot_{\Op} F_\p/\Op_\p\big) \ot_{\Op_\p}F_\p \lra T_p\big(E(L_\p)\ot_{\Op_\p} F_\p/\Op_\p\big) \ot_{\Op_\p}F_\p.
\end{equation}

We will be interested in the case when $L$ is an intermediate extension of a $\Zp$-extension $F_\infty$ of $F$. In the subsequent discussion, for a $F_\p[\Gal(F_\infty/F)]$-module $V$, we shall write $V[\Phi_n]$ to be the kernel of $F_\p$-endomorphism on $V$ given by multiplication by $\Phi_n$. We also sometimes call $V[\Phi_n]$ the $\Phi_n$-component (of $V$).

The next lemma elucidates the image of the map (\ref{f_L}) in the situation of our interest.

\bl \label{f_L lemma}
Let $F_\infty$ be a $\Zp$-extension of $F$ at which the prime $\p$ is totally ramified. Denote by $F_i$ the intermediate subextension of $F_\infty$ such that $|F_i:F|=p^i$. Let $n\geq 1$. Suppose that
\[ \rank_{\Op}E(F_n) - \rank_{\Op}E(F_{n-1}) >0.\]
Then we have
\[ \big(\im f_{F_n}\big)[\Phi_n] \cong F_\p[X]/\Phi_n,\]
where $f_{F_n}$ is defined as in $(\ref{f_L})$.
\el

\bpf
We follow the approach given in \cite[Lemma 4.2]{LeiZ}.
The logarithm map on the formal group of $E$ induces an isomorphism
 \[ T_p\big(E(F_{n,\p})\ot_{\Op_\p}F_\p/\Op_\p\big) \ot_{\Op_\p}F_\p \cong F_{n,\p} \cong F_\p[\Ga/\Ga_n] \]
of $F_\p[\Ga/\Ga_n]$-modules. From which, we have
\[ \Big(T_p\big(E(F_{n,\p})\ot_{\Op_\p}F_\p/\Op_\p\big) \ot_{\Op_\p}F_\p \Big) [\Phi_n] \cong F_\p[X]/\Phi_n. \]
Since $p$ is inert in $F/\Q$, the quadratic extension $F_\p/\Qp$ is unramified. Thus, the cyclotomic polynomial $\Phi_n$ remains irreducible over $F_\p[X]$ which in turn implies that $F_\p[X]/\Phi_n$ is an irreducible $\Ga/\Ga_n$-module. Therefore, it suffices to show that the restriction of the map $f_{F_n}$ on the $\Phi_n$-component is nonzero. By hypothesis, there is a nontorsion point $P\in E(F_n)-E(F_{n-1})$. Since the map $f_{F_n}$ respects the $\Ga/\Ga_n$-action, it sends the point $P$ to a nontorsion point in $E(F_{n,\p})-E(F_{n-1,\p})$ which yields the required conclusion. \epf

We are now ready to prove the first main theorem of the section.

\bt \label{mainCMcyc}
Suppose that the prime $p$ is inert in $F/\Q$. Let $F^\cyc$ be the cyclotomic $\Zp$-extension of $F$. Assume that $\Zhe(E/F_n)[p^\infty]$ is finite for every $n$, where $F_n$ is the subextension of $F^\cyc/F$ such that $|F_n:F|=p^n$. Then we have
\begin{equation}\label{mainCMcyceqn}
\mM(E/F^\cyc)^\vee\sim \bigoplus_{n\geq 0} (\La_{\Op_\p}/\Phi_n)^{\oplus \max\{0,f_n-1\}},
\end{equation}
where
\[ f_n = \left\{
           \begin{array}{ll}
            \displaystyle\frac{\rank_{\Op_F}E(F_n) - \rank_{\Op_F}E(F_{n-1}) }{\phi(p^n)} , & \mbox{if $n\geq 1$;} \\
            \\
            \rank_{\Op_F}E(F) , & \mbox{if $n=0$.}
           \end{array}
         \right.
\]
\et

\bpf
 It is well-known that $E(F^\cyc)$ is finitely generated over $\Z$ (see \cite{K, Ro}). Consequently, we have $f_n=0$ for $n\gg0$ and so the right hand side of $(\ref{mainCMcyceqn})$ is actually a finite sum. Furthermore, the finite generation allows us to apply an $\Op_{\p}$-analog of Proposition \ref{fineMWgeneral2} to obtain
\[ \mM(E/F^\cyc)^\vee\sim \bigoplus_{n\geq 0} (\La_{\Op_\p}/\Phi_j)^{\oplus s_n} \]
with the property that
\[T_p(\mM(E/F_n)) \cong \bigoplus_{j= 0}^n (\La_{\Op_{\p}}/\Phi_j)^{\oplus s_j}\]
for every $n$. It therefore remains to show that $s_n = \max\{0,f_n-1\}$ for every $n$. Now consider the following sequence
 \[ 0\lra \mM(E/F_n) \lra E(F_n)\ot\Qp/\Zp\lra E(F_{n,\p})\ot\Qp/\Zp. \]
 To relate the terms in this exact sequence to the map (\ref{f_L}), we first verify the following isomorphisms
 \[ E(F_n)\ot\Qp/\Zp\cong E(F_n)\ot_{\Op} F_\p/\Op_\p; \]
 \[ E(F_{n,\p})\ot\Qp/\Zp\cong E(F_{n,\p})\ot_{\Op_\p} F_\p/\Op_\p. \]
 Indeed, one can check easily that in the following diagram 
 \[   \entrymodifiers={!! <0pt, .8ex>+} \SelectTips{eu}{}\xymatrix{
    E(F_n)\ot_\Z\Zp \ar[d] \ar[r] & E(F_n)\ot\Qp   \ar[d] \ar[r] & E(F_n)\ot\Qp/\Zp\ar[d] \ar[r] & 0\\
       E(F_n)\ot_{\Op} \Op_\p \ar[r]^{} & E(F_n)\ot_{\Op_\p} F_\p \ar[r] & E(F_n)\ot_{\Op_\p} F_\p/\Op_\p \ar[r] & 0
     } \]
     that the leftmost vertical map is injective and the middle vertical map is bijective. Hence this yields the first claimed isomorphism. The second isomorphism can be established via a similar argument. In view of these identifications, upon applying the Tate module functor and tensoring with $F_\p$, we obtain an exact sequence
 \[ 0\lra T_p\big(\mM(E/F_n)\big)\ot_{\Op_\p} F_\p \lra T_p\big(E(F_n)\ot_{\Op} F_\p/\Op_\p\big)\ot_{\Op_\p} F_\p \stackrel{f_{F_n}}{\lra} T_p\big(E(F_{n,\p})\ot_{\Op_\p} F_\p/\Op_\p\big)\ot_{\Op_\p} F_\p, \]
 where $f_{F_n}$ is defined as in (\ref{f_L}).
 Taking the kernel of $\Phi_n$ and applying $\Op_\p$-analogs of Propositions \ref{fineMWgeneral2}, \ref{MWgrowth} and \ref{MWgrowthlocal}, we obtain an exact sequence
 \[ 0\lra \big(F_\p[X]/\Phi_n\big)^{\oplus s_n}\lra \big(F_\p[X]/\Phi_n\big)^{\oplus f_n}\stackrel{g_n}{\lra} F_\p[X]/\Phi_n,\]
 where the map $g_n$ is induced by the map $f_{F_n}$.
 If $f_n> 0$, then the map $g_n$ is surjective by Lemma \ref{f_L lemma} and so the above exact sequence is actually short exact. Therefore, we have $s_n=f_n-1$,  which therefore concludes the proof of the theorem.
\epf

\medskip \noindent
\textbf{Question 1.} Write $\La_{\Op_\p}=\Op_{\p}\ps{\Ga}$. Is the following formula
 \[\mathrm{char}_{\La_{\Op_\p}}Y(E/F^\mathrm{cyc}) =\bigoplus_{n\geq 0} \Phi_n^{\max\{0,f_n-1\}} \]
 valid? In view of Theorem \ref{mainCMcyc}, this is equivalent to asking whether the pseudo-isomorphism
 \[Y(E/F^\mathrm{ac})\sim \bigoplus_{n\geq 0} (\La_{\Op_\p}/\Phi_n)^{\oplus \max\{0,f_n-1\}}\] holds
or whether $\Zhe(E/F^{\cyc})$ is finite? Note that the question of the finiteness of this latter group has previously been posed by Wuthrich (see \cite{WuTS}; also see \cite[Question 2]{LimFineDoc}).

\medskip
We now consider the anticyclotomic situation.

\bt \label{mainCManticyc}
Suppose that the prime $p$ is inert in $F/\Q$. Let $F^\mathrm{ac}$ be the anticyclotomic $\Zp$-extension of $F$. Assume that $\Zhe(E/F_n)[p^\infty]$ is finite for every $n$, where $F_n$ is the subextension of $F^\mathrm{ac}/F$ such that $|F_n:F|=p^n$. Then we have
\begin{equation}\label{mainCManticyceqn}
Y_f(E/F^\mathrm{ac})\sim \bigoplus_{n\geq 0} (\La_{\Op_\p}/\Phi_n)^{\oplus \max\{0,f_n-1\}},
\end{equation}
where
\[ f_n = \left\{
           \begin{array}{ll}
            \displaystyle\frac{\rank_{\Op_F}E(F_n) - \rank_{\Op_F}E(F_{n-1}) }{\phi(p^n)} , & \mbox{if $n\geq 1$;} \\
            \\
            \rank_{\Op_F}E(F) , & \mbox{if $n=0$.}
           \end{array}
         \right.
\]
\et

\bpf
Note that in the settings of the theorem, one has $f_n = 0$ or $1$ for $n\gg 0$ accordingly to $n$ having parity equal to the functional equation of $L(E/\Q,s)$ or not (see \cite[Theorem 1.8 and discussion after formula (1.10)]{G01}). Therefore, the right hand side of (\ref{mainCManticyceqn}) is a finite sum and makes sense.
The proof is similar to that in Theorem \ref{mainCMcyc} except that in this context, $E(F^\mathrm{ac})$ is not finitely generated over $\Z$ and so Proposition \ref{fineMWgeneral2} cannot applied. But thankfully we do know that $Y(E/F^\ac)$ is a torsion $\La$-module by \cite[Proposition 3.3(i)]{AH2} and so we may appeal to Proposition \ref{fineMWgeneral} to obtain a result on $Y_f(E/F^\ac)$.
\epf

\br
In \cite{LeiQ}, Lei has proposed an analogue of \cite[Problem 0.7]{KurPo} in the context of a $p$-supersingular CM elliptic curve over the anticyclotomic $\Zp$-extension. In view of our result, Lei's formulation might be slightly incorrect. We therefore propose the following modification. \er

\medskip \noindent
\textbf{Question 2.} Write $\La_{\Op_\p}=\Op_{\p}\ps{\Ga}$. Is the following formula
 \[\mathrm{char}_{\La_{\Op_\p}}Y(E/F^\mathrm{ac}) =\bigoplus_{n\geq 0} \Phi_n^{\max\{0,f_n-1\}} \]
 valid? In view of Theorem \ref{mainCManticyc}, this is equivalent to asking whether the pseudo-isomorphism
 \[Y(E/F^\mathrm{ac})\sim \bigoplus_{n\geq 0} (\La_{\Op_\p}/\Phi_n)^{\oplus \max\{0,f_n-1\}}\] holds.

\br We briefly mention some interesting implications of this question. Recall that the signed Selmer group (see \cite{AH2, Kob, LeiQ} for definitions) is exactly non-cotorsion at the one whose sign parity agrees with that of the sign of the functional equation of $L(E/\Q,s)$ and cotorsion for the other (see \cite[Theorem 3.6]{AH2}). Therefore, if the assertion of Question 2 holds, it will follow from \cite[Theorem 3.5]{LeiQ} that the characteristic ideal of the cotorsion signed Selmer group consists precisely of product of cyclotomic polynomials. For more related discussion on this, we refer readers to \cite{LeiQ}.\er

\section{CM elliptic curve: good ordinary reduction case} \label{CMsplit section}

We continue to let $E$ denote an elliptic curve defined over $\Q$ with complex multiplication given by $\Op$, where $\Op$ is the ring of integers of an imaginary quadratic field $F$. Note that $F$ necessarily has class number $1$. In this section, \textbf{we shall assume that $E$ has good ordinary reduction at $p$, and so our prime $p$ splits in $F/\Q$}. Denote by $\p$ and $\overline{\p}$ the primes of $F$ above $p$. Since $F$ has class number 1, we have $\p = \pi\Op$ and $\overline{\p} = \overline{\pi}\Op$. We shall frequently identify $\Op_\p$ and $\Op_{\overline{\p}}$ with $\Zp$, and $F_\p$ and $F_{\overline{\p}}$ with $\Qp$ without any further mention.
We then define
\[ E[\p^\infty] = \cup_{n\geq 1} E[\pi^n] \quad\mbox{and}\quad E[\overline{\p}^\infty] = \cup_{n\geq 1} E[ \overline{\pi}^n]. \]
Note that $E[\p^\infty]$ and $E[\overline{\p}^\infty]$ are cofree $\Zp$-modules of corank $1$, and there is a natural decomposition
\[ \Ep= E[\p^\infty]\oplus E[\overline{\p}^\infty]\]
of $\Gal(\overline{\Q}/F)$-modules.
For a given $\Op$-module $M$, we shall write $T_\p(M) :=\plim_n M[\pi^n]$ and $T_{\overline{\p}}(M): =\plim_n M[\overline{\pi}^n]$. One can check easily that $T_pM = T_\p(M)\oplus T_{\overline{\p}}(M)$.


Let $\q\in\{\p, \overline{\p}\}$ and let $L$ be a finite extension of $F$. Following Coates \cite[P.\ 110, (6)]{C83}, we define the $\q$-primary fine Selmer group $R_\q(E/L)$ to be
\[ R_\q(E/L) =\ker\Big(H^1(G_S(L), E[\q^\infty]\Big)\lra \bigoplus_{w\in S(L)}H^1(L_w, E[\q^\infty])\Big).\]

\br
It is interesting to note that this $\q$-primary fine Selmer group has appeared in an 1983 article of Coates \cite[P.\ 110, (6)]{C83} in a different guise (also see the next remark). This was slightly more than 20 years before the paper of Coates-Sujatha \cite{CS05} which christened the fine Selmer group as we know today.
\er

We can also define the $\q$-primary Selmer group $\Sel_\q(A/L)$ by the exact sequence
\[0\lra \Sel_\q(A/L)\lra H^1(G_S(L),E[\q^\infty])\lra \bigoplus_{w\in S(L)} H^1(L_w, E)[\q^\infty].\]
Analogously to the usual Selmer group, the $\q$-primary Selmer group sits in the middle of the following short exact sequence
\[0\lra E(L)\ot_{\Op} F_\q/\Op_\q \lra \Sel_\q(E/L) \lra \sha(E/L)[\q^\infty] \lra 0.\]
Via a similar argument to that in \cite[Lemma 4.1]{LMu}, we have
\[ 0\lra R_\q(E/L) \lra \Sel_\q(E/L) \lra \bigoplus_{w\in S_\q(L)}E(L_w)\ot_{\Op_\q}F_\p/\Op_\q,\]
where $S_\q(L)$ is the set of primes of $L$ above $\q$.

\br In \cite{C83}, Coates has used the above exact sequence to define his $\q$-primary fine Selmer group (also compare with \cite{WuTS}). \er

We then define the fine $\q$-Mordell-Weil group by
\[ \mathcal{M}_\q(E/L) = \ker\Big(E(L)\ot_{\Op} F_\q/\Op_\q \lra \bigoplus_{w\in S_\q(L)} E(L_w)\ot_{\Op_\q}F_\p/\Op_\q \Big)\]
and the fine $\q$-Tate-Shafarevich group by
\[ \Zhe_\q(E/L) = \coker\Big( \mM_\q(E/L)\lra R_\q(E/L)\Big).\]
By an entirely similar discussion to that in Section 3, one can show that the $\q$-primary fine Selmer group sits in the middle of the following short exact sequence
\[0\lra \mM_\q(E/L) \lra R_\q(E/L) \lra \Zhe_\q(E/L) \lra 0\]
with $\Zhe_\q(E/L)$ injecting into $\sha(E/L)[\q^\infty]$.

 Suppose that now $L$ is a Galois extension of $F$ at which $\q$ is totally ramified in $L/F$.
Then the natural inclusion
$E(L) \hookrightarrow E(L_{\q})$ induces a $\Qp[\Gal(L/F)]$-map
\begin{equation}\label{f_Lsplit}
 f_L : T_\q\big(E(L)\ot_{\Op} F_\q/\Op_\q \big) \ot_{\Zp} \Qp \lra T_\q\big(E(L_\q)\ot_{\Op} F_\q/\Op_\q\big) \ot_{\Zp} \Qp.
\end{equation}
Note that here we have made use of the identification $\Op_\q\cong \Zp$ and $F_\q\cong \Qp$.

\bl \label{f_L lemma split}
Let $F_\infty$ be a $\Zp$-extension of $F$ at which the prime $\q$ is totally ramified. Denote by $F_i$ the intermediate subextension of $F_\infty$ such that $|F_i:F|=p^i$. Let $n\geq 1$. Suppose that
\[ \rank_{\Op}E(F_n) - \rank_{\Op}E(F_{n-1}) >0.\]
Then we have
\[ \big(\im f_{F_n}\big)[\Phi_n] \cong \Qp[X]/\Phi_n,\]
where $f_{F_n}$ is defined as in $(\ref{f_Lsplit})$.
\el

\bpf
This is proven similar to that in Lemma \ref{f_L lemma} (also see \cite[Lemma 4.2]{LeiZ}).
\epf

\bt \label{mainCMpsplit1}
Let $\q\in\{\p, \overline{\p}\}$. Denote by $F_{\q^\infty}$ the $\Zp$-extension of $F$ unramified outside $\q$. Assume that $\Zhe_\q(E/F_n)$ is finite for every $n$, where $F_n$ is the subextension of $F_{\q^\infty}/F$ such that $|F_n:F|=p^n$. Then there is a pseudo-isomorphism
\[\mM_\q(E/F_{\q^\infty})^\vee\sim \bigoplus_{n\geq 0} (\La/\Phi_j)^{\oplus \max\{0,f_n-1\}}\]
of $\Zp\ps{\Gal(F_{\q^\infty}/F)}$-modules,
where
\[ f_n = \left\{
           \begin{array}{ll}
            \displaystyle\frac{\rank_{\Op}E(F_n) - \rank_{\Op}E(F_{n-1}) }{\phi(p^n)} , & \mbox{if $n\geq 1;$} \\
            \\
            \rank_{\Op}E(F) , & \mbox{if $n=0$.}
           \end{array}
         \right.
\]
\et

\bpf
The approach of proof is essentially similar to those before. We shall therefore only explain the necessary modification and fine points. First of all, one can show that
\[ R_\q(E/F_n) \lra R_\q(E/F_{\q^\infty})^{\Ga_n}\]
has finite kernel and cokernel for every $n$ by an entirely similar argument to that in \cite[Theorem 3.3]{LimFineDoc}. A classical result of Coates (cf.\ \cite[Theorem 16]{C83}; also see \cite[Chap IV, Corollary 1.8]{DeSha}) tells us that $E(F_{\q^\infty})$ is finitely generated over $\Z$. Furthermore, the extension $F_{\q^\infty}/F$ is totally ramified at $\q$ (see \cite[Chap II, Proposition 1.9]{DeSha}) and so we may proceed as in Theorem \ref{mainCMcyc} applying Lemma \ref{f_L lemma split}.
\epf

In view of Theorem \ref{mainCMpsplit1}, we may pose the following question.

\medskip \noindent
\textbf{Question 3.} Is the following formula
\[ \mathrm{char}_{\Zp\ps\Ga}\Big(R_\q(E/F_{\q^\infty})^\vee\Big) = \prod_{n\geq 0} \Phi_n^{\max\{0,f_n-1\}} \]
valid? Note that by Theorem \ref{mainCMpsplit1}, this is also equivalent to asking whether
\[ R_\q(E/F_{\q^\infty})^\vee \sim \prod_{n\geq 0} (\La/\Phi_n)^{\oplus \max\{0,f_n-1\}}\]
or whether $\Zhe_{\q}(E/F_{\q^\infty})$ is finite (compare with \cite[Question 8.3]{WuTS}).

\smallskip

For the cyclotomic $\Zp$-extension, we have the following.

\bt \label{mainCMpsplitcyc}
Assume that $\Zhe(E/F_n)$ is finite for every $n$, where $F_n$ is the subextension of $F^\cyc/F$ such that $|F_n:F|=p^n$. Then we have
\[\mM(E/F^\cyc)^\vee\sim \bigoplus_{j\geq 0} (\La/\Phi_j)^{\oplus 2\max\{0,e_n-1\}},\]
where
\[ e_n = \left\{
           \begin{array}{ll}
            \displaystyle\frac{\rank_{\Z}E(F_n) - \rank_{\Z}E(F_{n-1}) }{\phi(p^n)} , & \mbox{if $n\geq 1$;} \\
            \\
            \rank_{\Z}E(F) , & \mbox{if $n=0$.}
           \end{array}
         \right.
\]
\et

\bpf
Observe that $\mM(E/F^\cyc) = \mM_{\p}(E/F^\cyc)\oplus \mM_{\overline{\p}}(E/F^\cyc)$.
Moreover, a theorem of Rubin \cite[Theorem 4.4]{Rub} asserts that $E(F^\cyc)$ is finitely generated over $\Z$. Therefore, we may apply a similar argument as before to obtain
 \[\mM_\q(E/F_{\q^\infty})^\vee\sim \bigoplus_{n\geq 0} (\La/\Phi_j)^{\oplus \max\{0,e_n-1\}}\]
 for $\q\in\{\p, \overline{\p}\}$. The conclusion of the theorem then follows up by patching up these conclusions.
\epf

The preceding theorem leads us to formulate the following question.

\medskip \noindent
\textbf{Question 4.} Is the following formula
\[ \mathrm{char}_{\Zp\ps\Ga}\Big(Y(E/F^{\cyc})\Big) = \prod_{n\geq 0} \Phi_n^{2\max\{0,e_n-1\}} \]
valid? In view of Theorem \ref{mainCMpsplitcyc}, this is equivalent to asking whether
 \[ Y(E/F^{\cyc}) \sim \prod_{n\geq 0} (\La/\Phi_n)^{\oplus 2\max\{0,e_n-1\}}\]
 or whether $\Zhe(E/F^{\cyc})$ is finite.

\medskip
We finally come to the anticyclotomic $\Zp$-extension.

\bt \label{mainCMpsplitac}
Assume that $\Zhe(E/F_n)$ is finite for every $n$, where $F_n$ is the subextension of $F^\ac/F$ such that $|F_n:F|=p^n$. Define
\[ e_n = \left\{
           \begin{array}{ll}
            \displaystyle\frac{\rank_{\Z}E(F_n) - \rank_{\Z}E(F_{n-1}) }{\phi(p^n)} , & \mbox{if $n\geq 1$;} \\
            \\
            \rank_{\Z}E(F) , & \mbox{if $n=0$.}
           \end{array}
         \right.
\]
\begin{itemize}
  \item[$(1)$] If the root number of $E/\Q$ is $+1$, then we have
\[ \mM(E/F^\ac)^\vee\sim \bigoplus_{j\geq 0} (\La/\Phi_j)^{\oplus 2\max\{0,e_n-1\}}.\]
  \item[$(2)$] If the root number of $E/\Q$ is $-1$, then we have
\[ Y_f(E/F^\ac)\sim \bigoplus_{j\geq 0} (\La/\Phi_j)^{\oplus 2\max\{0, e_n-1\}}.\]
\end{itemize}
\et

\bpf
If the root number of $E/\Q$ is $+1$, then a classical theorem of Greenberg
\cite[Theorem 3]{G83} asserts that the group $E(F^\ac)$ is finitely generated and we may proceed as before appealing to Proposition \ref{fineMWgeneral2}. If the root number of $E/\Q$ is $-1$, we at least know that $Y(E/F^\ac)$ is torsion over $\La$ (cf. \cite[Proposition 2.4.4]{AH}) and so Proposition \ref{fineMWgeneral} applies.
\epf

In view of our result, we therefore pose the following question.

\medskip \noindent
\textbf{Question 5.} Is the following formula
\[ \mathrm{char}_{\Zp\ps\Ga}\Big(Y(E/F^{\ac})\Big) = \prod_{n\geq 0} \Phi_n^{2\max\{0,e_n-1\} } \]
valid?

\section{Growth of fine Mordell-Weil groups under generalized Heegner hypothesis} \label{BDP section}

In this section, $E$ will always be an elliptic curve defined over $\Q$ which has good reduction at $p\geq 5$. In other words, $E$ either has good ordinary reduction or good supersingular reduction at $p$. Let $F$ be an imaginary quadratic field at which the prime $p$ splits completely and does not divide the class number of $F$. We shall always assume our data $(E,F)$ satisfies the generalized Heegner hypothesis which we now recall.
The conductor $N$ of the elliptic curve $E$ will have a factorization $N= N^+N^-$ such that every prime dividing $N^+$ splits in $F/\Q$, whereas every prime dividing $N^-$ is inert in $F/\Q$. We further assume that $N^-$ is a squarefree product of an even number of primes. Note that in this context, the group $E(F^\ac)$ is not finitely generated over $\Z$ (see \cite{CW, JSW}). In the event that $E$ has good ordinary reduction, we assume further that $E(\mathbb{F}_p)$ has no $p$-torsion. (Note that the latter is automatically true had the reduction type of $E$ at $p$ is supersingular; see \cite[Proposition 8.7]{Kob}.)
Finally, we always assume that the natural group homomorphism $\Gal(\overline{\Q}/\Q)\lra \mathrm{Aut}_{\Zp}(T_pE)$ induced by the Galois action on $E$ is surjective. All the assumptions in this paragraph will be in force throughout this section.

The goal of this section is to attempt to study the variation of the fine Mordell-Weil group in the intermediate subextensions of $F^\ac/F$, and hence the structure of $Y_f(E/F^\ac)$, under the above setting. Throughout this section, we shall let $F_n$ denote the subextension of $F^\ac/F$ such that $|F_n:F|=p^n$. We also define
\[ e_n = \left\{
           \begin{array}{ll}
            \displaystyle\frac{\rank_{\Z}E(F_n) - \rank_{\Z}E(F_{n-1}) }{\phi(p^n)} , & \mbox{if $n\geq 1;$} \\
            \\
            \rank_{\Z}E(F) , & \mbox{if $n=0$.}
           \end{array}
         \right.
\]

Let $\p$ and $\overline{\p}$ be the primes of $F$ above $p$. Since the prime $p$ is assumed to not divide the class number of $F$, both $\p$ and $\overline{\p}$ are totally ramified in $F^\ac/F$. In this situation, we have two localization maps
\[ \psi_{n,\p}:E(F_n)\ot\Qp/\Zp \lra E(F_{n,\p})\ot\Qp/\Zp,   \]
\[\psi_{n,\overline{\p}}:E(F_n)\ot\Qp/\Zp \lra E(F_{n,\overline{\p}})\ot\Qp/\Zp.\]

\bl \label{onto local map}
Assume that $\sha(E/F_n)[p^\infty]$ is finite for every $n$, where $F_n$ is the subextension of $F^\ac/F$ such that $|F_n:F|=p^n$. Then the maps $\psi_{n,\p}$ and $\psi_{n,\overline{\p}}$ are surjective for every $n$.
\el

\bpf
When $E$ has good supersingular reduction at $p$, this is established in \cite[Lemma 5.3]{LLM}. We shall therefore consider the case when $E$ has good ordinary reduction at $p$. By the argument in Lemma \ref{f_L lemma}, it suffices to show that $e_n>0$ for every $n$. But the latter follows from an application of \cite[Remark 2.1.6(2)]{Lee}.
\epf

We now proceed to prove an analog of our results for the BDP-Selmer group in the sense of \cite{JSW} using the techniques develop in this paper. We keep the standing assumptions as in the first paragraph of this section.
We recall two variants of Selmer groups which were introduced in \cite{JSW}, although we shall adhere to the notation used in \cite{LLM}. The first of which is defined as follow
\begin{multline*}
  \Sel^{\mathrm{Gr}, 0}(E/F_n) =\ker\Bigg(H^1(G_S(F_n),\Ep)\lra \\
  H^1(F_{n,\p},E)[p^\infty]\oplus H^1(F_{n,\overline{\p}},\Ep)\oplus \Big(\bigoplus_{w\in S(F_n), w\nmid p} H^1(F_{n,w},\Ep)\Big)\Bigg).
\end{multline*}

The next is defined by
\[\Sel^{\mathrm{BDP}}(E/F_n) =\ker\Big(H^1(G_S(F_n),\Ep)\lra H^1(F_{n,\overline{\p}},\Ep)\oplus\big(\bigoplus_{w\in S(F_n), w\nmid p} H^1(F_{n,w},\Ep)\Big).\]
Set $\Sel^{\mathrm{Gr}, 0}(E/F^\ac) = \ilim_n \Sel^{\mathrm{Gr}, 0}(E/F_n)$ and $\Sel^{\mathrm{BDP}}(E/F^\ac) = \ilim_n \Sel^{\mathrm{BDP}}(E/F_n)$. Then write $X^{\mathrm{Gr}, 0}(E/F^\ac)$ for the Pontryagin dual of $\Sel^{\mathrm{Gr}, 0}(E/F^\ac)$ and
\[ X_f^{\mathrm{Gr}, 0}(E/F^\ac) = \plim_n X^{\mathrm{Gr}, 0}(E/F_n)/X^{\mathrm{Gr}, 0}(E/F_n)[p^\infty].\]
The modules $X^{\mathrm{BDP}}(E/F^\ac)$  and $X_f^{\mathrm{BDP}}(E/F^\ac)$ are defined similarly.
The connection between these Selmer groups is summarized in the next lemma.

\bl \label{Sel compare}
For each $n$, we have the following exact sequences
\[ 0\lra \ker \psi_{n,\overline{\p}} \lra \Sel^{\mathrm{Gr}, 0}(E/F_n) \lra \sha(E/F_n)[p^\infty]\lra 0,\]
\[ 0\lra \Sel^{\mathrm{Gr}, 0}(E/F_n)\lra \Sel^{\mathrm{BDP}}(E/F_n)  \stackrel{h_n}{\lra} H^1(F_{n,\p},E)[p^\infty],\]
where $\mathrm{im} ~h_n$ is finite.
\el

\bpf
When the elliptic curve has good supersingular reduction, the first exact sequence is \cite[Proposition 5.5]{LLM}, and the second exact sequence and the finiteness of $\mathrm{im} ~h_n$ can be found in the proof of \cite[Proposition 5.6]{LLM}. An examination of the proof reveals that the global duality argument there carries over to the good ordinary case.
\epf

We record certain useful consequences of the preceding lemma in the form of the following proposition.

\bp \label{control SelGr}
One has $X_f^{\mathrm{Gr}, 0}(E/F^\ac) = X_f^{\mathrm{BDP}}(E/F^\ac)$.
 Furthermore, the restriction map
\[ \Sel^{\mathrm{Gr}, 0}(E/F_n)  \lra \Sel^{\mathrm{Gr}, 0}(E/F^\ac)^{\Ga_n}\]
is injective with a finite cokernel.
\ep

\bpf
Taking Pontryagin dual of the second exact sequence in Lemma \ref{Sel compare}, we obtain
\[ 0 \lra C_n\lra  X^{\mathrm{BDP}}(E/F_n) \lra X^{\mathrm{Gr}, 0}(E/F_n) \lra 0, \]
for some finite $C_n$. Since $C_n$ is finite, we have the following exact sequence
\[ 0 \lra C_n\lra  X^{\mathrm{BDP}}(E/F_n)[p^\infty] \lra X^{\mathrm{Gr}, 0}(E/F_n)[p^\infty] \lra 0 \]
by \cite[Lemma 2.1]{LimFineDoc}. Combining these two exact sequences, we obtain an identification
\[ X^{\mathrm{BDP}}(E/F_n)\big/X^{\mathrm{BDP}}(E/F_n)[p^\infty] \cong X^{\mathrm{Gr}, 0}(E/F_n)\big/X^{\mathrm{Gr}, 0}(E/F_n)[p^\infty]. \]
Taking inverse limit yields the first equality of the proposition.

For the final assertion of the proposition, consider the following commutative diagram
\[   \entrymodifiers={!! <0pt, .8ex>+} \SelectTips{eu}{}\xymatrix{
    0 \ar[r]^{} & \Sel^{\mathrm{Gr}, 0}(E/F_n) \ar[d] \ar[r] & \Sel^{\mathrm{BDP}}(E/F_n)
    \ar[d] \ar[r] & \mathrm{im}~h_n  \ar[d] \ar[r] & 0\\
    0 \ar[r]^{} & \Sel^{\mathrm{Gr}, 0}(E/F^\ac)^{\Ga_n} \ar[r]^{} & \Sel^{\mathrm{BDP}}(E/F^\ac)^{\Ga_n} \ar[r] & \big(\ilim_n \mathrm{im}~h_n\big)^{\Ga_n} &
     } \]
     with exact rows. By \cite[Theorem 4.1]{LLM}, the middle map is injective with a finite cokernel. Since $\im~h_n$ is also finite by Proposition \ref{Sel compare}, an application of the snake lemma then shows that the leftmost vertical map is injective with finite cokernel. This completes the proof of the proposition.
\epf

\br
In general, it would seem that the cokernels in the preceding proposition need not be bounded independently of $n$  (see \cite[Proposition 5.6 and Corollary 6.12]{LLM}). Nevertheless, the finiteness suffices for our purposes (see hypothesis (b) of our Corollary \ref{tech module}).
\er

\br
Although we do not require this fact, we like to mention that in the event that $E$ has good supersingular reduction at $p$, one even has
$\Sel^{\mathrm{Gr}, 0}(E/F^\ac) =\Sel^{\mathrm{BDP}}(E/F^\ac)$. Indeed, this assertion is a consequence of  $\ilim_n H^1(F_{n,\p},E)[p^\infty] =0$ and the second exact sequence in Proposition \ref{Sel compare}.
\er

We now come to the following theorem.

\bt \label{BDPac}
Assume that $\sha(E/F_n)[p^\infty]$ is finite for every $n$, where $F_n$ is the subextension of $F^\ac/F$ such that $|F_n:F|=p^n$. 
Then we have
\[  X_f^{\mathrm{BDP}}(E/F^\ac) \sim \bigoplus_{n\geq 0} (\La/\Phi_n)^{\oplus (e_n-1)}.\]
\et

\bpf
By \cite[Corollary 3.5]{Cas} and \cite[Theorem 4.4]{LLM}, $\Sel^{\mathrm{BDP}}(E/F^\ac)$ is cotorsion over $\La$. In view of Lemma \ref{Sel compare}, Proposition \ref{control SelGr} and the finiteness of $\sha(E/F_n)[p^\infty]$, we may apply Corollary \ref{tech module} to conclude that
\[  X_f^{\mathrm{BDP}}(E/F^\ac)=X_f^{\mathrm{Gr}, 0}(E/F^\ac) \sim \bigoplus_{n\geq 0} (\La/\Phi_j)^{\oplus t_n}\]
with the property that at each level $n$, one has
\[ T_p(\ker\psi_{n,\overline{\p}}) \cong \bigoplus_{j=0}^n (\La/\Phi_j)^{\oplus t_j}.\]
It remains to show that $t_n=e_n-1$ for every $n$. By virtue of Lemma \ref{onto local map}, there is the following short exact sequence
\[ 0 \lra \ker\psi_{n,\overline{\p}} \lra E(F_n)\ot\Qp/\Zp \lra E(F_{n,\p})\ot\Qp/\Zp \lra 0. \]
Now, comparing $\Zp$-coranks of the terms in this said exact sequence, we obtain
\[ \sum_{j=0}^n e_j\deg (\Phi_j) = \sum_{j=0}^n t_j\deg (\Phi_j) + \sum_{j=0}^n \deg (\Phi_j) = \sum_{j=0}^n (t_j+1)\deg (\Phi_j)  \]
for every $n$. The desired equality follows from this and an inductive argument.
\epf

\br \label{BDPremark}
Recall that there is a main conjecture relating the BDP-Selmer group and the so-called BDP $p$-adic $L$-function \cite{BDP} (also see \cite{Cas,  CLW, CW, Kob-O, LZ} for recent progress in this main conjecture, where this list is far from being exhaustive). Therefore, the preceding theorem gives a link between the precise growth of Mordell-Weil groups and the BDP $p$-adic $L$-function. (See discussion in subsection below)
\er

We now come to our result on the fine Selmer groups.

\bt \label{BDPfine}
Assume that $\sha(E/F_n)[p^\infty]$ is finite for every $n$, where $F_n$ is the subextension of $F^\ac/F$ such that $|F_n:F|=p^n$.
Then we have
\[  Y_f(E/F^\ac) \sim \bigoplus_{n\geq 0} (\La/\Phi_n)^{\oplus s_n}\]
such that
\[ T_p (\mM(E/F_n)) \cong \bigoplus_{j = 0}^n (\La/\Phi_j)^{\oplus s_j} \]
and that
\[ \max\{0,e_n-2\} \leq s_n \leq \max\{0,e_n-1\}\]
for every $n$.
\et

\bpf
As seen in the proof of Theorem \ref{BDPac}, $\Sel^{\mathrm{BDP}}(E/F^\ac)$ is cotorsion over $\La$. Since $R(E/F^\ac)$ is contained in $\Sel^{\mathrm{BDP}}(E/F^\ac)$, it is also  cotorsion over $\La$. Therefore, we may apply Proposition \ref{fineMWgeneral2} to conclude that
\[  Y_f(E/F^\ac) \sim \bigoplus_{n\geq 0} (\La/\Phi_n)^{\oplus s_n}\]
such that for each $n$, one has
\[ T_p(\mM(E/F_n)) \cong \bigoplus_{j=0}^n (\La/\Phi_j)^{\oplus s_j}.\]

On the other hand, it follows from definition that there is the following exact sequence
\[ 0\lra \mM(E/F_n) \lra \ker\psi_{n,\overline{\p}} \lra E(F_{n,\p})\ot\Qp/\Zp. \]
Upon taking Tate-module functor, followed by tensoring with $\Qp$ and considering the $\Phi_n$-component, we obtain
\[ 0 \lra \big(\Qp[X]/\Phi_n\big) ^{\oplus s_n} \lra \big(\Qp[X]/\Phi_n\big) ^{\oplus (e_n-1)}\lra \Qp[X]/\Phi_n.\]
(Note that for the middle term here, we have made use of the fact that
\[ T_p(\ker\psi_{n,\overline{\p}}) \cong \bigoplus_{j=0}^n (\La/\Phi_j)^{\oplus (e_j-1)}\]
which was derived in the proof of Theorem \ref{BDPac}.) From which, we may conclude that
\[ \max\{0,e_n-2\} \leq s_n \leq \max\{0,e_n-1\}\]
as required.
\epf

Inspired by the above theorem, we pose the following question.

\medskip \noindent
\textbf{Question 6.} Can $\mathrm{char}_{\Zp\ps\Ga}\big(Y(E/F^{\mathrm{ac}})\big)$ be entirely described by cyclotomic polynomial? More precisely, is it possible that the following formula
\[ \mathrm{char}_{\Zp\ps\Ga}\Big(Y(E/F^{\mathrm{ac}})\Big) = \prod_{n\geq 0} \Phi_n^{s_n} \]
holds for some $s_n$ satisfying
\[ \max\{0,e_n-2\} \leq s_n \leq \max\{0,e_n-1\}?\]

\subsection{Appendix: More on Remark \ref{BDPremark}}

 The goal of this short subsection is to say a bit on Remark \ref{BDPremark}. We continue to retain the settings of the section. Now, let $\mathfrak{L}_\mathrm{BDP}$ denote the BDP $p$-adic $L$-function attached to $E$ over $F$ (in the sense of \cite{BDP}). Note that $\mathfrak{L}_\mathrm{BDP}$  is an element in
$\La^{\mathrm{ur}}$, where $\La^{\mathrm{ur}} := \La \otimes \Zp^{\mathrm{ur}}$ and $\Zp^{\mathrm{ur}}$ is the $p$-adic completion of the ring of integers
of the maximal unramified extension of $\Qp$. We now state the following conjecture.

\begin{conjecture} \label{BDPconj}
$\mathfrak{L}_\mathrm{BDP}^2 \in \mathrm{char}_{\La}\Big(X^{\mathrm{BDP}}(E/F^\mathrm{ac})\Big)\La^{\mathrm{ur}}$.
  \end{conjecture}

This conjecture is essentially one-side inclusion of the BDP main conjecture, and it has been verified
in many cases (see \cite{Kob-O, LZ}).

\bc \label{BDPacCor}
Retain the assumptions in Theorem \ref{BDPac}. Suppose further that Conjecture \ref{BDPconj}
is valid. Then we have an inequality
\[\ord_{\Phi_n}\Big(\mathfrak{L}_\mathrm{BDP}\Big)\geq \frac{e_n-1}{2}.\]
\ec

\bpf
Indeed, it follows from Theorem \ref{BDPac} and the assumption of the corollary that
\[2~\ord_{\Phi_n}\Big(\mathfrak{L}_\mathrm{BDP}\Big)\geq \ord_{\Phi_n}\left(\mathrm{char}_{\La}\Big(X^{\mathrm{BDP}}(E/F^\mathrm{ac})\Big)\right) \geq   \ord_{\Phi_n}\left(\mathrm{char}_{\La}\Big(X_f^{\mathrm{BDP}}(E/F^\mathrm{ac})\Big)\right) =e_n-1,\]
where the second inequality follows from $X_f^{\mathrm{BDP}}(E/F^\mathrm{ac})$ being a quotient of $X^{\mathrm{BDP}}(E/F^\mathrm{ac})$.
\epf

\end{document}